\documentclass[graybox]{svmult}

\usepackage{type1cm}        
\usepackage{makeidx}         
\usepackage{graphicx}        
\usepackage{multicol}        
\usepackage[bottom]{footmisc}

\usepackage{newtxtext}       %
\usepackage{newtxmath}       

\renewtheorem{example}{Example}

\def\eqref#1{(\ref{#1})}
\def\dsp{\displaystyle}
\def\Frac#1#2{\frac
{
 {\raise.6ex
 \hbox{$\displaystyle#1$}}
}
{
 {\lower.6ex
 \hbox{$\displaystyle#2$}}
 }
}

\numberwithin{equation}{section}

\makeatother

\def\dsp{\displaystyle}
\def\Frac#1#2{\frac
{
 {\raise.6ex
 \hbox{$\displaystyle#1$}}
}
{
 {\lower.6ex
 \hbox{$\displaystyle#2$}}
 }
}

\def\sign{{\rm sign}}

\def\CHFs#1#2#3{
{}_1F_1\left({a};{c};{z}\right)
}

\def\FG#1#2#3#4{
{}_2F_1\left(
\begin{array}{c}
\begin{array}{c}\hskip-10pt#1,#2\end{array}\\
\begin{array}{c}\hskip-10pt #3\end{array}
\end{array}
\hskip-8pt;\,#4
\right)}

\def\erfc{{\rm erfc}}

\def\bigO{{\cal O}}

\def\tfrac#1#2{{{\lower.6ex
\hbox{$\scriptstyle#1$}}\over 
{\raise.7ex
\hbox{$\scriptstyle#2$}}}}

\def\RR{\mathbb R}

\def\sign{{\rm sign}}

\def\erfc{{{\rm erfc}}}

\def\erfc{{\rm erfc}}

\def\tfrac#1#2{{{\lower.6ex
\hbox{$\scriptstyle#1$}}\over 
{\raise.7ex
\hbox{$\scriptstyle#2$}}}}

\makeindex             

\begin{document}

\title*{A new asymptotic representation and inversion method for the Student's $t$ distribution}

\author{A. Gil, J. Segura and N. M. Temme}
\institute{A. Gil \at Departamento de Matem\'atica Aplicada y CC. de la Computaci\'on.
ETSI Caminos. Universidad de Cantabria. 39005-Santander, Spain. \email{amparo.gil@unican.es}
\and J. Segura \at Departamento de Matem\'aticas, Estad\'{\i}stica y 
        Computaci\'on, Universidad de Cantabria, 39005 Santander, Spain. \email{javier.segura@unican.es}
\and N.M. Temme \at IAA, 1825 BD 25, Alkmaar, The Netherlands 
(Former address: CWI, Science Park 123, 1098 XG Amsterdam,  The Netherlands). \email{nicot@cwi.nl}}
%

%
%
\maketitle

\abstract{Some special functions are particularly relevant in applied probability 
and statistics. For example, the incomplete beta function is the 
cumulative central beta distribution. 
 In this paper, we consider the inversion
of the central Student's-$t$ distribution which is a particular case of the central beta distribution.
The inversion of this distribution functions is useful in hypothesis testing as well as for generating random samples distributed according to the corresponding probability density function.
A new asymptotic representation in terms of
the complementary error function, will be one of the important ingredients in our analysis. As we will show, this asymptotic representation is also useful in the computation of the distribution function. 
We illustrate the performance of all the  obtained approximations  with numerical examples.
}

\section{Introduction}\label{sec:int}

There is a very close relationship between some special functions and some of the most popular distribution functions in statistics.
For example, the incomplete beta function \cite[\S8.17]{Paris:2010:INC} is the central beta distribution. Particular cases include other well-known distributions such as the geometric, binomial, negative binomial or the central Student's-$t$ distribution. Therefore, standard methods for the computation and inversion 
of special functions \cite{Gil:2007:NSF} can also be applied to evaluate and invert distribution functions.  The problem of inversion appears, for example,
when computing percentage points of the distribution functions; also, it is closely related to the generation of random variates from a continuous probability density function needed, for example, in Monte Carlo or quasi-Monte Carlo methods.

We have considered the central beta distribution in a previous publication \cite{Gil:2017:IBE}.
In this paper, we focus on the computation and inversion of the central Student's-$t$ distribution, which has multiple applications in science and engineering (for an application in physics, see for example \cite{Rover:2011:RSD}). As we mentioned before, this distribution function is a particular case of the central beta 
distribution but which requires particular analysis.
An asymptotic representation in terms of
the complementary error function, will play a key role in our analysis.
The performance of the approximations obtained will be illustrated with numerical examples.

There is a vast literature on Student's $t$ distribution. An extensive overview can be found at \cite{Johnson:1995:CUD}; see also the references 
contained therein.
For an historical account and the origin of this distribution, we refer to \cite{Zabell:2008:OAA}. For  a generalization from the  viewpoint of special functions, see \cite{Koepf:2006:GSD}.

Some useful expressions for the analysis of the central Student's-t distribution are:

\begin{description}

\item{a) Probability density function:}

\begin{equation}\label{eq:intro01}
f_n(t)=\frac{\Gamma\left(\frac12n+\frac12\right)}{\sqrt{n\pi}\,\Gamma\left(\frac12n\right)}\left(1+\frac{t^2}{n}\right)^{-\frac12n-\frac12}=
\frac{1}{\sqrt{n}\,B\left(\frac12,\frac12n\right)}\left(1+\frac{t^2}{n}\right)^{-\frac12n-\frac12},
\end{equation}
where $ t\in{\RR}$ and $B(p,q)$ is the Beta integral; $n>0$, not necessarily an integer.

\item{b) Cumulative distribution function:}

\begin{equation}\label{eq:intro02}
F_n(x)=\int_{-\infty}^x f_n(t)\,dt,\quad x\in\RR.
\end{equation}

\item{c) Incomplete beta function:}

\begin{equation}\label{eq:intro03}
I_x(p,q)=\frac{1}{B(p,q)}\int_0^x t^{p-1}(1-t)^{q-1}\,dt, \quad B(p,q)=\frac{\Gamma(p)\Gamma(q)}{\Gamma(p+q)}.
\end{equation}
From the integral representation, we have:
\begin{equation}\label{eq:intro04}
I_x(p,q)=1-I_{1-x}(q,p).
\end{equation}

\item{d) Cumulative distribution function in terms of the incomplete beta function:}

\renewcommand{\arraystretch}{1.5}
\begin{equation}\label{eq:intro05}
F_n(x)=
\left\{
\begin{array}{ll}
\frac12+\frac12I_{\frac{x^2}{n+x^2}}\left(\frac12,\frac12n\right)=1-\frac12I_{\frac{n}{n+x^2}}\left(\frac12n,\frac12\right), \quad & {\rm if\ } x\ge0,\\
\frac12-\frac12I_{\frac{x^2}{n+x^2}}\left(\frac12,\frac12n\right)=\frac12I_{\frac{n}{n+x^2}}\left(\frac12n,\frac12\right), \quad & {\rm if\ } x\le0.\\
\end{array}
\right.
\renewcommand{\arraystretch}{1.0}
\end{equation}

\item{e) Cumulative distribution function in terms of the Gauss hypergeometric functions:}

\begin{equation}\label{eq:intro06}
F_n(x)=
\left\{
\begin{array}{ll}
\dsp{\tfrac12+\frac{x}{\sqrt{n}\,B\left(\frac12,\frac12n\right)}}\FG{\frac12}{\frac12n+\frac12}{\frac32}{-\frac{x^2}{n}},
 \quad & {\rm if\ } x\in\RR,\\
1-\dsp{\frac{(1-y)^{\frac12n}y^{-\frac12}}{nB\left(\frac12,\frac12n\right)}\FG{1}{\frac12}{\frac12n+1}{-\frac{n}{x^2}},}
 \quad & {\rm if\ }  x\ge0,\\
\dsp{\frac{(1-y)^{\frac12n}y^{-\frac12}}{nB\left(\frac12,\frac12n\right)}\FG{1}{\frac12}{\frac12n+1}{-\frac{n}{x^2}},}
 \quad & {\rm if\ }  x\le0,\\
\end{array}
\right.
\renewcommand{\arraystretch}{1.0}
\end{equation}
where $y=x^2/(n+x^2)$. The first formula in  \eqref{eq:intro06} is given in \cite{Amos:1964:RBD}, the other ones follow from well-known relations between the incomplete beta function and the hypergeometric function; see \cite[\S8.17(ii)]{Paris:2010:INC}.

\end{description}

\section{Asymptotic expansion of the Student's $t$ cumulative distribution function}\label{sec:as}

The second and third representation in \eqref{eq:intro06} can be used for large values of $n$ by using the standard power series of the hypergeometric functions. It is not necessarily that $x^2>n$, but a condition $x^2/n =\bigO(1)$ is needed. To obtain a large-$n$ asymptotic representation, whether or not $x$ is large, we use a method that we have used for other cumulative distribution functions; see  \cite[Chapter~36]{Temme:2015:AMI}.

We use in \eqref{eq:intro01}  the substitution $t=s\sqrt{n}$. This gives
\begin{equation}\label{eq:as01}
\begin{array}{@{}r@{\;}c@{\;}l@{}}
F_n(x)&=&\dsp{\frac{1}{\sqrt{n}\,B\left(\frac12,\frac12n\right)}\int_{-\infty}^{x}\left(1+\frac{t^2}{n}\right)^{-\frac12n-\frac12}\,dt}\\[8pt]
&=& \dsp{\frac{1}{B\left(\frac12,\frac12n\right)}
\int_{-\infty}^{x/\sqrt{n}}
\left(1+s^2\right)^{-\frac12n-\frac12}\,ds.}
\end{array}
\end{equation}
We write $u^2=\ln(1+s^2)$, with the condition $\sign(u)=\sign(s)$, and obtain
\begin{equation}\label{eq:as02}
F_n(x)=\frac{1}{B\left(\frac12,\frac12n\right)}\int_{-\infty}^\xi e^{-\frac12n u^2}g(u)\,du,
\end{equation}
where
\begin{equation}\label{eq:as03}
g(u)=\sqrt{\frac{u^2}{1-e^{-u^2}}}, \quad \xi^2=\ln\left(1+\frac{x^2}{n}\right),\quad \sign(\xi)=\sign(x).
\end{equation}

Using the method of  \cite[\S36.1]{Temme:2015:AMI} we find that $F_n(x)$ can be written in the form
\begin{equation}\label{eq:as04}
F_n(x)=\tfrac12\erfc\left(-\xi\sqrt{n/2}\right)-\frac{e^{-\frac12n\xi^2}}{\sqrt{2\pi n}\,\beta(n)}B_n(\xi),
\quad \beta(n)=\sqrt{\frac{n}{2\pi}}\,B\left(\tfrac12,\tfrac12n\right),
\end{equation}
where we have introduced the complementary error function
\begin{equation}\label{eq:as05}
\erfc\,z=\frac{2}{\sqrt{\pi}}\int_z^\infty e^{-t^2}\,dt.
\end{equation}

The function $B_n(\xi)$ has the asymptotic expansion
\begin{equation}\label{eq:as06}
B_n(\xi)\sim\sum_{k=0}^\infty \frac{C_k(\xi)}{n^k},\quad n\to\infty,\quad \xi\in\RR,
\end{equation}
where the coefficients follow from the recursive scheme
\begin{equation}\label{eq:as07}
C_k(\xi)=\frac{g_k(\xi)-g_k(0)}{\xi},\quad g_{k+1}(u)=\frac{d}{du}\frac{g_k(u)-g_k(0)}{u},
\end{equation}
with $g_0(u)=g(u)$ defined in \eqref{eq:as03}. The first coefficients are
\begin{equation}\label{eq:as08}
\begin{array}{@{}r@{\;}c@{\;}l@{}}
C_0(\xi)&=&\dsp{\frac{g(\xi)-1}{\xi}  },\\[8pt]
C_1(\xi)&=& \dsp{ -\frac{4g(\xi)^3-4g(\xi)\xi^2+\xi^2-4}{4\xi^3}},\\[8pt]
C_2(\xi)&=& \dsp{\frac{96g(\xi)^5-128g(\xi)^3\xi^2+32g(\xi)\xi^4-\xi^4+8\xi^2-96}{32\xi^5}}.\end{array}
\end{equation}

The function $\beta(n)$ has the expansion
\begin{equation}\label{eq:as09}
\beta(n)\sim\sum_{k=0}^\infty \frac{D_k}{n^k},\quad n\to\infty, \quad D_k=g_k(0).
\end{equation}
The first coefficients are
\begin{equation}\label{eq:as10}
D_0=1,\quad D_1=\tfrac14,\quad D_2=\tfrac{1}{32},\quad D_3=-\tfrac{5}{128}.
\end{equation}
These coefficients can be expressed in terms of the coefficients $a_k$ of $\dsp{g(u)=\sum_{k=0}^\infty a_k u^k}$. We have
(see \cite[Remark~36.2]{Temme:2015:AMI})
\begin{equation}\label{eq:as11}
D_k=\left(\tfrac12\right)^k 2^k a_{2k},\quad k=0,1,2,\ldots\,.
\end{equation}

\begin{figure}
\begin{center}
\includegraphics[width=13cm]{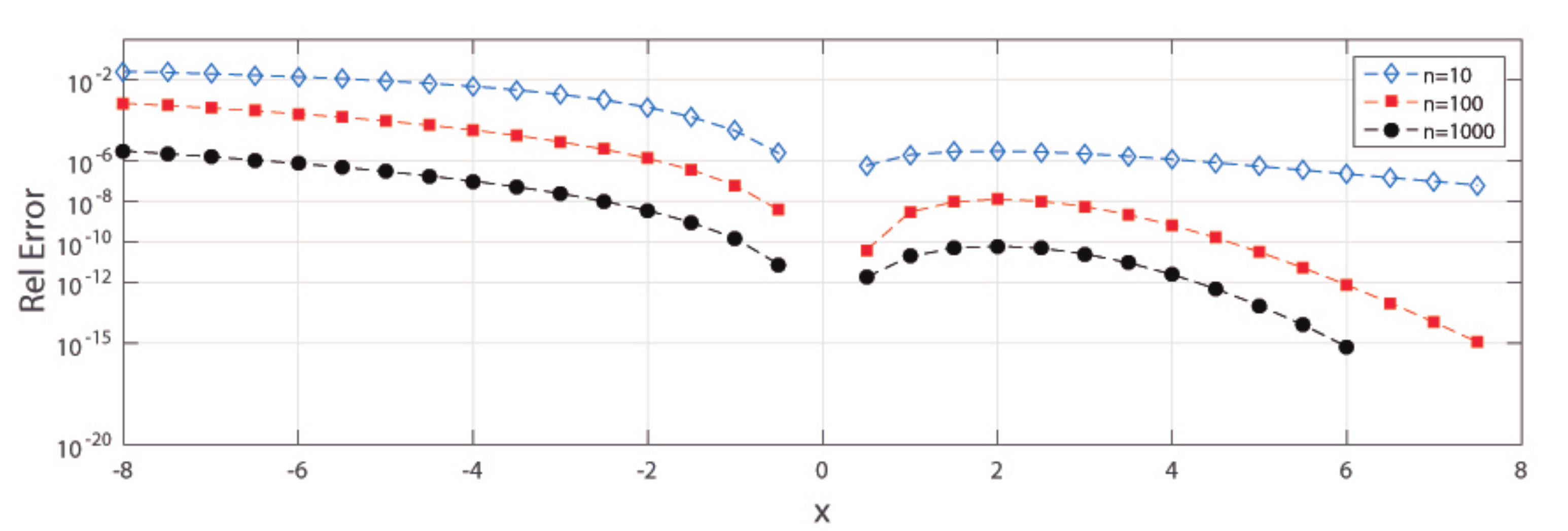}
\caption{
\label{fig:fig00} Relative errors obtained when computing the Student's $t$ cumulative 
distribution function using the expansion  \eqref{eq:as04}. }
\end{center}
\end{figure}

Examples of the performance of the expansion \eqref{eq:as04} for three values of $n$ ($n=10,\,100,\,1000$) are shown in Figure \ref{fig:fig00}.
Five $C_k$coefficients in the series \eqref{eq:as06} have been considered in the computations. Relative errors in comparison to the values
of the distribution function given in \eqref{eq:intro06} computed with Maple.

\section{Inversion of the Student's $t$ cumulative distribution function}\label{sec:invstud}
The inversion problem is: find $x$ that satisfies the equation
\begin{equation}\label{eq:ninv01}
F_n(x)=p,\quad 0 < p < 1. 
\end{equation}
We consider three different approaches, the first one is for small values of $p-\frac12$, which will give small values of $x$. The second method is for small values of $p$, which gives large negative values $x$. Thirdly we use the uniform approximation of \S\ref{sec:as}, which will be valid for a large range of $x$, including the  values near $x=0$.

\subsection{Inversion for small values of $p-\frac12$}\label{sec:sq}
We use the first representation in \eqref{eq:intro06}, and write the inversion problem in the form
\begin{equation}\label{eq:ninv02}
x\sum_{k=0}^\infty c_kx^{2k}=q,\quad q=\left(p-\tfrac12\right)\sqrt{n}\,B\left(\tfrac12,\tfrac12n\right),\end{equation}
and the $c_k$ follow from the coefficients of the hypergeometric function:
\begin{equation}\label{eq:ninv03}
c_k=\frac{\left(\frac12\right)_k\left(\frac12n+\frac12\right)_k}{k!\,\left(\frac32\right)_k}\frac{(-1)^k}{n^k}.\end{equation}
The solution of the equation in \eqref{eq:ninv01} has the expansion $\dsp{x=q\sum_{k=0}^\infty x_k q^{2k}}$, and the first coefficients are
\begin{equation}\label{eq:ninv04}
\begin{array}{@{}r@{\;}c@{\;}l@{}}
x_0&=&1,\quad \dsp{x_1=\frac{n+1}{6n},\quad  x_2=\frac{(n+1)(7n+1)}{120n^2},}\\[8pt]
 x_3&=&\dsp{\frac{(n+1)(127n^2+8n+1)}{5040n^3},
\quad
x_4=\frac{(n + 1)(4369n^3 - 537n^2 + 135n + 1)}{362880n^4}.} 
\end{array}
\end{equation}

\begin{figure}
\begin{center}
\includegraphics[width=13cm]{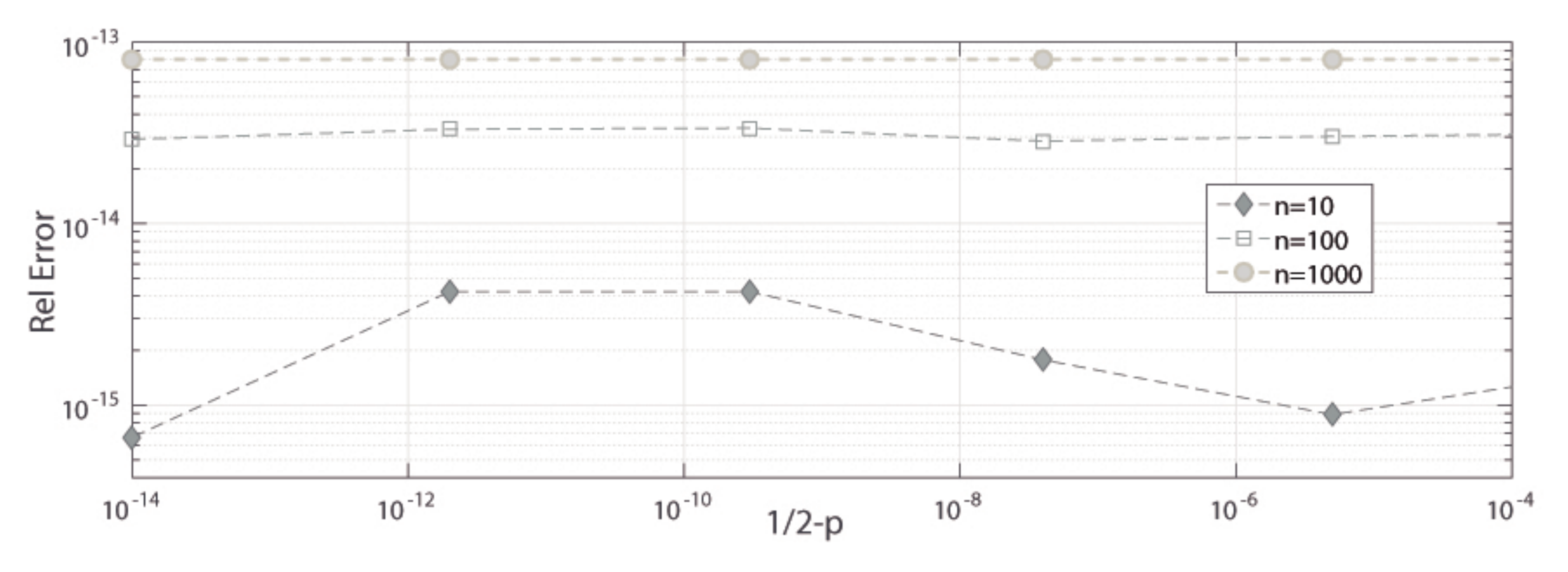}
\caption{
\label{fig:fig01} Performance of the expansion  for small values of $p-\frac12$ for three values of $n$. The five coefficients given in \eqref{eq:ninv04}
for the expansion have been considered in the computations. }
\end{center}
\end{figure}

We see that the shown coefficients are bounded for large values of $n$. 
Also, $q=\left(p-\frac12\right)\sqrt{2\pi}\,\beta(n)$ (see \eqref{eq:as04}) with an expansion of $\beta(n)$ given in \eqref{eq:as09}. This shows that the inversion considered here for small values of $\left(p-\frac12\right)$  is rather well conditioned for large values of $n$.
In \eqref{fig:fig01} we show examples of the performance of the expansion (using the four terms  given in \eqref{eq:ninv04}) 
for three values of $n$. The values of $p$ considered in the calculations
are $p=\frac12-\Delta$, for $\Delta=10^{-14}$, $2\times 10^{-12}$, $3\times 10^{-10}$, $4\times 10^{-8}$, $5\times 10^{-6}$, $6\times 10^{-4}$. The results obtained with the expansion have been compared with the Matlab function for the inversion of the central Student's-t distribution (function {\bf tinv}). As can be seen, a relative error better than $10^{-13}$ is obtained in all cases.

\begin{remark}\label{rem:01}
By using the definition in \eqref{eq:intro02}, it is easily verified that  the function $F_n(x)$ becomes for $n=1$
\begin{equation}\label{eq:ninv05}
F_1(x)=\frac12+\frac{1}{\pi}\arctan\,x,
\end{equation}
and for the inversion problem $F_n(x)=p$ we have $q=\left(p-\frac12\right)\pi$. The  equation to be inverted  becomes
$\arctan x=q$, with solution $x=\tan q$. The shown coefficients in \eqref{eq:ninv04} correspond with those of the expansion $x=q+\frac13q^3+\frac{2}{15}q^5+\dots$. This expansion converges for $0<p<1$.
\end{remark}

\subsection{High-order iteration}\label{sec:iter}
It should be mentioned that it is also possible to obtain numerical approximations for the inversion problem for small values of $p-\frac12$ using the fixed point iterations giving sharp error bounds for the central beta distribution described in \cite{Gil:2017:IBE}. For example, iterating the fixed point iterations $y=g(y)$ or $y=h(y)$  where 

\begin{equation}\label{eq:ninvb1}
\begin{array}{l}
g(y)=\left(2(p-\frac12)B(\frac12, \frac12n)(\frac12-(\frac12+\frac12n)y)\right)^2  (1-y)^{-n},\\[8pt]
h(y)=225\Frac{\left(p-\frac12\right)^2B\left(\frac12,\frac12n\right)^2 (1-y)^{-n}    }{\left(y^2(n^2+4n+3)+5y(n+1)+15\right)^2},
\end{array}
\end{equation}
starting from $y=0$, approximations to the values of $x$ in \eqref{eq:ninv01} are obtained with $x=\displaystyle\sqrt{\Frac{ny}{1-y}}$. As an example, using 2 iterations of the fixed point iteration $y=h(y)$ for $n=10$, $p=0.5+10^{-5}$ we find $x=0.00002569978035$  with a relative error $5\times 10^{-12}$.
An explicit expression for the $x$ value obtained in the second iteration is

\begin{equation}\label{eq:ninvb2}
x=120\left(\Frac{nZ(1-(1/4)Z)^{-n}}{Z^4A_4+Z^3A_3+Z^2A_2+ZA_1+57600}\right)^{1/2},
\end{equation}

where $Z=\left(2(p-\frac12)B(\frac12, \frac12n)\right)^2$ and

\begin{equation}\label{eq:ninvb3}
\begin{array}{l}
A_1=-14400(1-(1/4)Z)^{-n}+9600n+9600,\\
A_2=880n^2+2720n+1840,\\
A_3=40n^3+200n^2+280n+120,\\
A_4=n^4+8n^3+22n^2+24n+9.
\end{array}
\end{equation}

\begin{figure}
\begin{center}
\includegraphics[width=13cm]{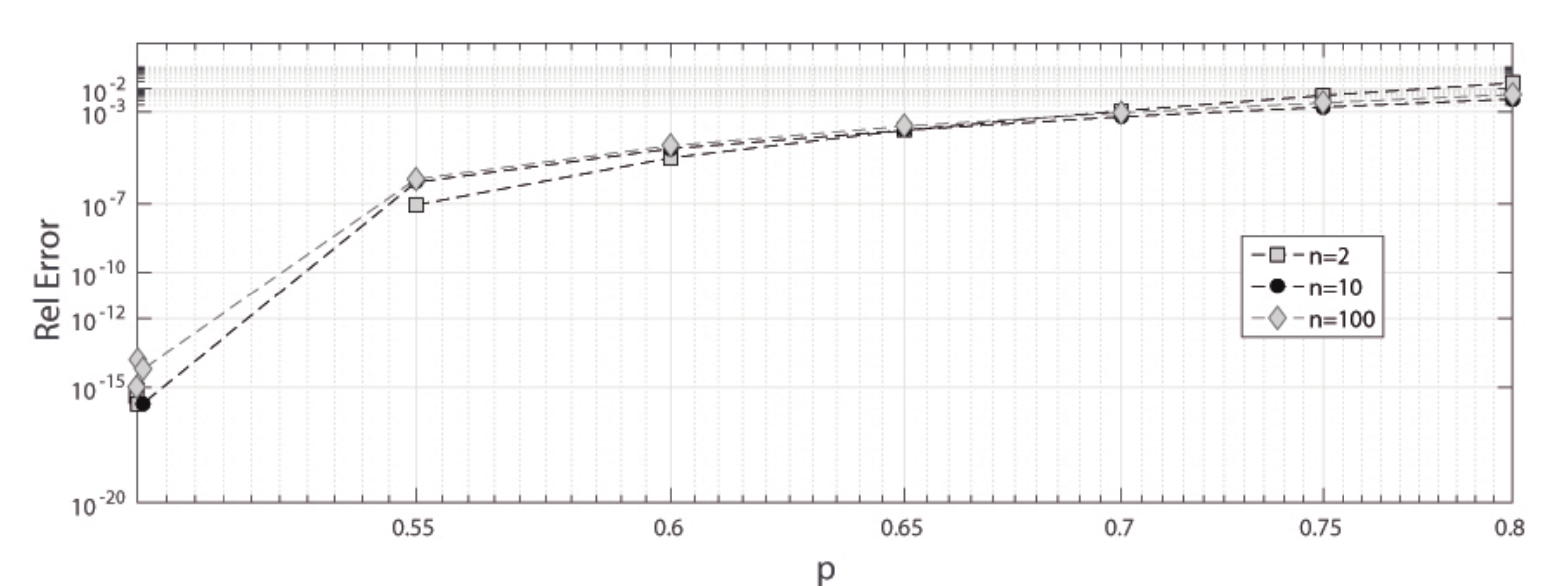}
\caption{
\label{fig:fig01b} Relative errors obtained when using \eqref{eq:ninvb2} to approximate the solution to the inversion problem
\eqref{eq:ninv01} for $n=2,\,10,\,100$.}
\end{center}
\end{figure}

This approximation can be also used for not so small values of $p-\frac12$, as can be seen in Figure \ref{fig:fig01b}.     

\subsection{Inversion for small values of $p$}\label{sec:sp}
For this case we use\footnote{The approach of this section is similar as one of the inversion methods considered in \cite{Hill:1970:StQ}.} the representation in the third line of \eqref{eq:intro06} and introduce
\begin{equation}\label{eq:ninv06}
\eta=\frac{n}{x^2},\quad \delta= \left(pnB\left(\tfrac12,\tfrac12n\right)\right)^{\frac{2}{n}}.
\end{equation}
Then we can write the equation $F_n(x)=p$ in the form
\begin{equation}\label{eq:ninv07}
\eta(1+\eta)^{\frac{1}{n}-1}S(\eta)^{\frac{2}{n}}=\delta,
\end{equation}
where $S(\eta)$ is the standard power series of the hypergeometric function in the third line of   \eqref{eq:intro06} .

We see that for small values of $\delta$ the wanted variable $\eta$ behaves as $\eta\sim\delta$, and substituting the expansion
\begin{equation}\label{eq:ninv08}
\eta=\sum_{k=1}^\infty\eta_k \delta^k,  
\end{equation}
we find  the following first few coefficients
\begin{equation}\label{eq:ninv09}
\begin{array}{l}
\eta_1=1,\quad \eta_2=\Frac{n+1}{n+2},\quad \eta_3=\Frac{(n+1)\left(2n^2 +9n + 6\right)}{2(n+2)^2(n+4)},\\[8pt]
\eta_4= \Frac{3n^6+35n^5+134n^4+328n^3+1174n^2+2100n+1152}{3n(n+2)^3(n+4)(n+6)}.\\
\end{array}
\end{equation}
When we have computed $\eta$, $x$ follows from \eqref{eq:ninv06}: $x=-\sqrt{n/\eta}$.

Again, we see that the coefficients satisfy $\eta_k=\bigO(1)$ for large values of $n$, which also happens for all coefficients that we have evaluated. The first 10 coefficients satisfy $\eta_k=1+\bigO(1/n)$ as $n\to\infty$. When $n$ is large, the only problem is that $\delta$ tends to~1, which for convergence of the series in \eqref{eq:ninv08} is a bad condition.

For example, when $p=10^{-8}$ and $n=10$, we have $\delta\doteq 0.038193186$. With 5 terms in the series in \eqref{eq:ninv08} we find $\eta\doteq0.039576861$, giving $x\doteq -15.8956879$. Then, $F_n(x)=9.9999981\times10^{-9}$, with relative error $1.92\times10^{-7}$.
For $n=25$ we find $\delta\doteq0.281$, $x\doteq-8.0759$ and a relative error $0.011$.
Other examples of the performance of the series  for small values of $p$
are given  in Table~\ref{table1}, where we show  the values and relative errors obtained for $n=10$ and few values of $p$.
The four coefficients given in \eqref{eq:ninv09} have been considered in the calculations.
 Since the Matlab function {\bf tinv} seems to fail for very small values of $p$, the tests have been performed comparing with
Maple.

\renewcommand{\arraystretch}{1.25}
\begin{table}
$$
\begin{array}{ccc}
p & x & \mbox{Rel. Error} \\
\hline
1\times10^{-50} & -256452.5718769479 & 6.5 \times 10^{-15} \\
2\times10^{-40} & -23927.87084268530 & 3.2 \times 10^{-15} \\
 3\times10^{-30} &  -2297.706518629116 &1.0 \times 10^{-13} \\
 4\times 10^{-20} & -223.234400503956 &1.1 \times 10^{-09} \\
 5\times10^{-10} & -21.62201646469524 &1.2 \times 10^{-05} \\
\end{array}
$$
\caption{
\label{table1} Inversion values $x$ obtained (for $n=10$) using the series in \eqref{eq:ninv08} for small values of $p$.  
The coefficients given in \eqref{eq:ninv09}  have been considered 
in the calculations. Relative errors in comparison to Maple, are also shown in the table. }
\end{table}
\renewcommand{\arraystretch}{1.0}

\begin{remark}\label{rem:02}
When we take $n=1$ and use \eqref{eq:ninv05}, the inversion problem becomes
\begin{equation}\label{eq:ninv10}
x=-\cot(p\pi),\quad \eta=\tan^2(p\pi),\quad \delta=p^2\pi^2, \quad 0<p<1.
\end{equation}
The expansion in \eqref{eq:ninv08} becomes for $p<\frac12$ and $\delta <\frac14\pi^2$
\begin{equation}\label{eq:ninv11}
\tan^2(p\pi)= \sum_{k=1}^\infty\eta_k \delta^k,\quad \eta_1=1,\quad \eta_2=\tfrac23,\quad \eta_3= \tfrac{17}{45}.
\end{equation}
These values of $\eta_k$ correspond with those given in \eqref{eq:ninv09} for $n=1$.
\end{remark}

\subsection{Inversion by using the uniform expansion}\label{sec:asuni}
We use the representation given in \eqref{eq:as04} and first try to find $\xi$, then $x$ follows from the relation in \eqref{eq:as03}.
We assume that $n$ is large and use  the asymptotic  method as described in \cite[\S42.1]{Temme:2015:AMI} and in our papers 
\cite{Gil:2017:IBE}, \cite{Gil:2019:NCB},  \cite{Gil:2020:IBD}.

Let $\xi_0$ satisfy the equation 
\begin{equation}\label{eq:ninv12}
\tfrac12\erfc\left(-\xi_0\sqrt{n/2}\right)=p.
\end{equation} 
Then we assume for $\xi$ the expansion
\begin{equation}\label{eq:ninv13}
\xi\sim\xi_0+\frac{\xi_1}{n}+\frac{\xi_2}{n^2}+\frac{\xi_3}{n^3}+\ldots,
\end{equation}
where $\xi_k$ have to be determined. When we have this approximation $\xi$ we compute $x$ from \eqref{eq:as03}.

From \eqref{eq:ninv01},  \eqref{eq:ninv12} and  \eqref{eq:as02} we find
\begin{equation}\label{eq:ninv14}
\frac{dp}{d\xi_0}=\sqrt{\frac{n}{2\pi}} e^{-\frac12n\xi_0^2},\quad
\frac{dp}{d\xi}=\sqrt{\frac{n}{2\pi}} \frac{g(\xi)}{\beta(n)}e^{-\frac12n\xi^2},
\end{equation}
where $\beta(n)$ is defined in \eqref{eq:as04}.

Dividing the two derivatives, we find
\begin{equation}\label{eq:ninv15}
g(\xi)\frac{d\xi}{d\xi_0}=\beta(n)e^{\frac12n(\xi^2-\xi_0^2)},
\end{equation}
We substitute the expansion given  in \eqref{eq:ninv13}, use $\beta(n)=1+\bigO(1/n)$, and obtain, considering equal large-order terms of $n$, the next term in the expansion:
\begin{equation}\label{eq:ninv16}
g(\xi_0)=e^{\xi_0\xi_1} \quad \Longrightarrow \quad \xi_1=\frac{1}{\xi_0}\ln g(\xi_0).
\end{equation}
Because 
\begin{equation}\label{eq:ninv17}
g(u)=1+\tfrac14 u^2+\tfrac{1}{96}u^4+\bigO\left(u^6\right),
\end{equation}
it follows that $\xi_1$ is well defined when $\xi_0$ tends to zero  (that is, when $p\sim\frac12$).

We can find higher-order terms $\xi_j, j\ge 2,$ of the expansion in \eqref{eq:ninv13} using more coefficients in the asymptotic expansion of $\beta(n)$ (see  \eqref{eq:as10}). Also, we need the expansion
\begin{equation}\label{eq:ninv18}
g(\xi)=g(\xi_0)+(\xi-\xi_0)g^{\prime}(\xi_0)+\tfrac12(\xi-\xi_0)^2g^{\prime\prime}(\xi_0)+\ldots .
\end{equation}
By using \eqref{eq:ninv15} and algebraic manipulations we find a few other coefficients:
\begin{equation}\label{eq:ninv19}
\begin{array}{@{}r@{\;}c@{\;}l@{}}
\xi_2&=&-\Bigl(2g\xi \xi_1^2  + 4\left(g-\xi g^{\prime}\right)\xi_1 + \xi  g - 4g^\prime\Bigr)/(4\xi^2g),\\[8pt]
\xi_3&=&
\Bigl(2\xi^2g^ 2\xi_1^3 +\left(2\xi^3gg^{\prime\prime}  - 2\xi^3{g^\prime}^2 -6\xi^2gg^\prime+  8\xi g^2\right)\xi_1^2+\\
&&\left(12g + \xi^2g  - 16\xi g^\prime+ 4\xi^2g^{\prime\prime}\right) g\xi_1+\\ 
&& \xi g ^2 + 4\xi gg^{\prime\prime} + 2\xi {g^\prime}^2- \xi^2gg^\prime - 12gg^\prime   \Bigr)\Bigl/  (4\xi^4g^2),
\end{array}
\end{equation}
where $\xi=\xi_0$ and $g$, $g^\prime$ and $g^{\prime\prime}$ are evaluated at $\xi_0$.

For small values of $\xi_0$ (that is, when $p\sim\frac12$), we need expansions. We have

\begin{equation}\label{eq:ninv20}
\begin{array}{@{}r@{\;}c@{\;}l@{}}
\xi_1&=&\dsp{\ \tfrac{1}{4}\xi_0-\tfrac{1}{48}\xi_0^3+\tfrac{1}{5760}\xi_0^7-\tfrac{1}{362880}\xi_0^{11}+\tfrac{1}{19353600}\xi_0^{15}+
\ldots,}\\[8pt]
\xi_2&=&\dsp{\tfrac{1}{32}\xi_0-\tfrac{5}{192}\xi_0^3+\tfrac{7}{2560}\xi_0^5+\tfrac{1}{2560}\xi_0^7-\tfrac{407}{5806080}\xi_0^9-\tfrac{13}{1451520}\xi_0^{11}+\ldots,}\\[8pt]
\xi_3&=&\dsp{ -\tfrac{5}{128}\xi_0-\tfrac{11}{1536}\xi_0^3+\tfrac{63}{10240}\xi_0^5-\tfrac{823}{2580480}\xi_0^7-\tfrac{5291}{23224320}\xi_0^
9+\ldots,}\\[8pt]
\xi_4&=&\dsp{ -\tfrac{21}{2048}\xi_0+\tfrac{37}{2048}\xi_0^3+\tfrac{179}{81920}\xi_0^5-\tfrac{22711}{10321920}\xi_0^7+\ldots,}\\[8pt]
\xi_5&=&\dsp{ \tfrac{399}{8192}\xi_0+\tfrac{219}{32768}\xi_0^3-\tfrac{3679}{327680}\xi_0^5+\ldots,}\\[8pt]
\xi_6&=&\dsp{\tfrac{869}{65536}\xi_0-\tfrac{6877}{131072}\xi_0^3+\ldots\ .}
\end{array}
\end{equation}

\begin{example}\label{exam:01}
We summarise the algorithmic steps for the inversion method. We take $n=10$,  $p=0.44$ and two terms in the expansion in \eqref{eq:ninv13}.
\begin{enumerate}
\item
Compute $\xi_0$ from equation \eqref{eq:ninv12}. We have $\xi_0\doteq-0.047746$.
\item
Compute $\xi_1$ from equation \eqref{eq:ninv16} by using $g(u)$ defined in \eqref{eq:as03}: $\xi_1\doteq -0.011933$.
\item
Compute $\xi$ by using \eqref{eq:ninv13} with the information now available: $\xi\sim \xi_0+\xi_1/n\doteq  -0.048934$.
\item
Compute $x$ from equation \eqref{eq:as03}. Because $\xi<0$,  $x$ should be negative: $x \doteq-0.1548354$.
\item
Verification: compute $F_n(x)$ by using the first line in \eqref{eq:intro06}: $F_{10}(x)\doteq 0.4400158$. Relative error: $0.000036$.
\end{enumerate}
\end{example}

\begin{figure}
\begin{center}
\includegraphics[width=13cm]{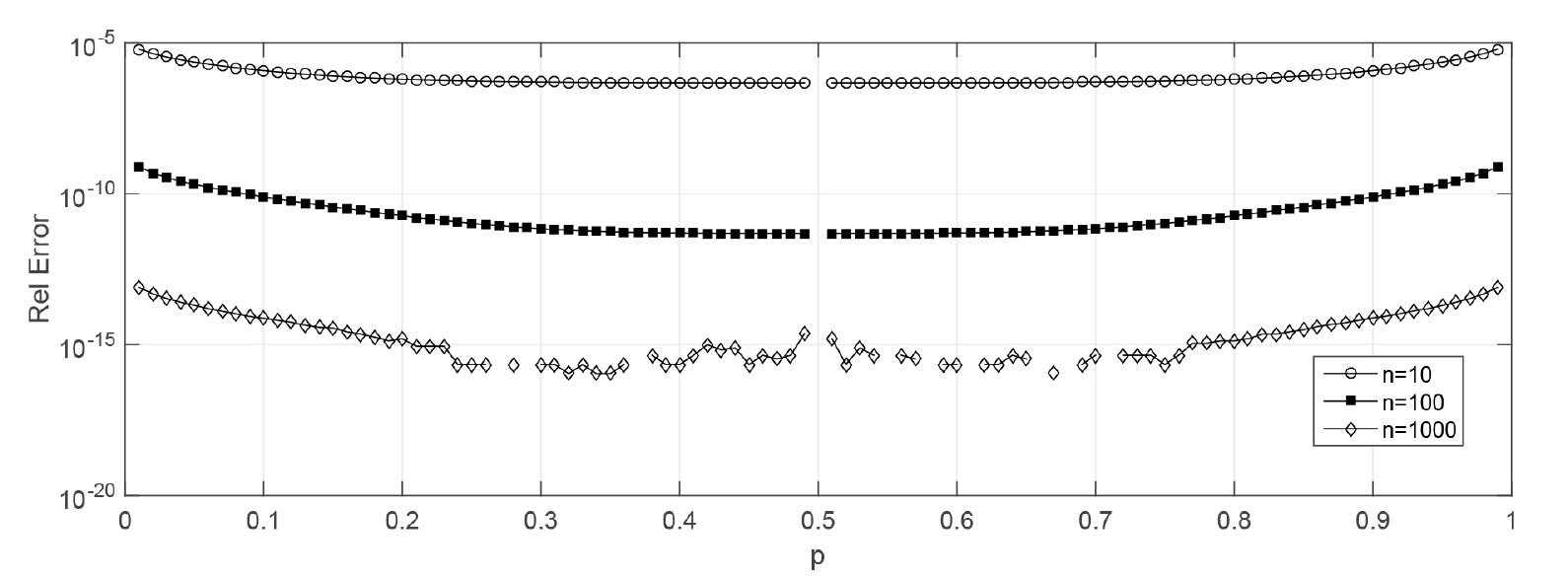}
\caption{
\label{fig:fig02} Performance of the asymptotic inversion method  using four terms of the expansion for small values of $\xi_0$ given 
in \eqref{eq:ninv20}.  Three  values of $n$ are considered in the computations. }
\end{center}
\end{figure}

A test of the performance of the asymptotic inversion method for three values of $n$ using the expansion for small values of $\xi_0$ given in \eqref{eq:ninv20},
 can be seen in Figure \ref{fig:fig02}. The terms of the expansion $\xi_1,\xi_2, \xi_3,\xi_4$ given in \eqref{eq:ninv20} 
have been considered in the computations.  For computing the inverse of the complementary error function needed to
compute $\xi_0$ in \label{eq:ninv12} we use the function {\bf inverfc} given in \cite{Gil:2015:GCH}.
The relative errors obtained in comparison to the Matlab function {\bf tinv} are shown in the figure.

\section{Concluding remarks}

We have presented approximations for the inversion problem \eqref{eq:ninv01} of the central Student-$t$ distributions. To obtain the approximations, different methods have been considered depending on the values of $p$. In particular, one of the key elements in our analysis was
 the use of an asymptotic representation of the distribution function in terms of the complementary error function. Numerical tests have shown that the approximations obtained are, in all cases, accurate. Also, they are easy to compute, which is
an important advantage when using the inverse to generate random variates  distributed according to central Student-$t$ probability density function.

\bibliographystyle{plain}
\bibliography{biblio}

\begin{acknowledgement}
We acknowledge financial support from Ministerio de Ciencia e Innovaci\'on, Spain, 
project PGC2018-098279-B-I00 (MCIU/AEI/FEDER, UE). 
NMT thanks CWI, Amsterdam, for scientific support.
\end{acknowledgement}

\end{document}